\documentclass[12pt,twoside]{article}    
\usepackage{amsfonts,amsmath,latexsym}

\newtheorem{thm}{Theorem}[section]

\newtheorem{lem}[thm]{Lemma}

\newtheorem{defn}[thm]{Definition}

\numberwithin{equation}{section}
\newenvironment{proof}[1][Proof]{\textbf{#1.} }{\ \rule{0.5em}{0.5em}}

\makeatletter \@addtoreset{equation}{section} \makeatother
\begin{document}


\pagestyle{myheadings}

\markboth{\hfill {\small AbdulRahman Al-Hussein } \hfill}{\hfill
{\small Maximum principle for optimal control of stochastic partial differential equations } \hfill }


\thispagestyle{plain}


\begin{center}
{\large \textbf{Maximum principle for optimal control of stochastic partial differential equations}$^{*}$\footnotetext{$^{*}$
This work is supported by the Science College Research Center at Qassim University, project no. SR-D-011-724. }} \\
\vspace{0.7cm} {\large AbdulRahman Al-Hussein }
\\
\vspace{0.2cm} {\footnotesize
{\it Department of Mathematics, College of Science, Qassim University, \\
 P.O.Box 6644, Buraydah 51452, Saudi Arabia \\ {\emph E-mail:} alhusseinqu@hotmail.com}}
\end{center}

\begin{abstract}
We shall consider a stochastic maximum principle of optimal control for a
control problem associated with a stochastic partial
differential equations of the following type:
\begin{eqnarray*}
 \left\{ \begin{array}{ll}
              d x(t) = ( A (t) x(t) + a (t , u (t) ) x(t) + b(t , u(t)) ) d t
              \\ \hspace{3cm} + \, [ \big< \sigma (t , u(t)) , x(t) \big>_{K} + \, g (t , u(t)) ] d M(t) , \\
             \; x(0) = x_0 \in K ,
         \end{array}
 \right.
\end{eqnarray*}
with some given predictable mappings $a , b , \sigma , g$ and a continuous martingale $M$ taking its values in a Hilbert space $K ,$
while $u(\cdot)$ represents a control. The equation is also driven by
a random unbounded linear operator $A(t ,w), \; t \in [0,T ], $ on $K .$

We shall derive necessary conditions of optimality for this control problem without a convexity assumption on
the control domain, where $u(\cdot)$ lives, and also when this control variable
is allowed to enter in the martingale part of the equation.
\end{abstract}



\section{Introduction}\label{sec:int}
Consider the following stochastic partial differential equation (SPDE for short):
\begin{eqnarray}\label{intr:spde}
 \left\{ \begin{array}{ll}
              d x(t) = ( A (t) x(t) + a (t , u (t) ) x(t) + b(t , u(t)) ) d t
              \\ \hspace{3cm} + \, [ \big< \sigma (t , u(t)) , x(t) \big>_{K} + \, g (t , u(t)) ] d M(t) , \, 0 \leq t \leq T, \\
             \; x(0) = x_0 \in K ,
         \end{array}
 \right.
\end{eqnarray}
where $A (t) , t \in [0 , T], $ is a random unbounded closed linear
operator on a separable Hilbert space $K .$ The noise is modelled by a continuous martingale $M$ in $K$ and $a , b , \sigma$ and $g$ are suitable predictable bounded mappings while $u (\cdot )$ is a control. This equation will be studied over a Gelfand triple $( V , K , V' ) .$  That is $V$ is
a separable Hilbert space embedded continuously and densely in $K .$
More precisely, given a bounded measurable mapping $\ell : [ 0 , T] \times \mathcal{O} \rightarrow K$ and a fixed element $G$ of $K ,$ we shall be interested here in minimizing the cost functional:
\begin{equation*}
J(u (\cdot ) ) = \mathbb{E} \; \big{[} \; \int_0^T \big< \ell (t , u (t) ) , x^{u
(\cdot ) } (t) \big>_{K} \, dt + \big< G , x^{u (\cdot ) } (T) \big>_{K} \; \big{]} ,
\end{equation*}
over the set of admissible controls. We will approach this by using the adjoint equation of the SPDE~(\ref{intr:spde}), which is a backward stochastic partial differential equation (BSPDE) driven by an infinite dimensional martingale, and derive in particular a stochastic maximum principle for this optimal control problem. Such BSPDEs (or even BSDEs) have their importance shown in applications in control theory like \cite{[Alh-AMO-10]} and in some financial applications as in \cite{[Marie09]}. For more applications we refer the reader to Bally et al. \cite{[Bal-Pard05]}, Imkeller et al. \cite{[Imk10]} and \cite{[Fuh-Tess]}.

It is known that a Wiener filtration is usually required to deal with BSPDEs that arise as adjoint equations of controlled SPDEs. This is indeed a restriction insisted on for example in \cite{[Y-Z]} and \cite{[Zhou93]}. {\O}ksendal et al. in \cite{[Oks-Z05]} and some other recent works have now considered the adjoint equation of a controlled BSPDE with a filtration generated by a Wiener process and a Poisson random measure. In our work here we can consider an arbitrary continuous filtration thanks to a result established in \cite{[Alh-Stoc09]} giving existence and uniqueness of solutions to BSPDEs driven by martingales. In this respect we refer the reader also to Imkeller et al. \cite{[Imk10]}, where a filtration is being taken which is similar to the one used here.
The reader can also see \cite{[Alh-AMO-10]}, \cite{[Tud89]}, \cite{[Kot84]}, \cite{[Gyo-Kr81]}, \cite{[Gyo83]}, \cite{[Kry-Roz07]}, \cite{[Gre-Tu95]}, \cite{[Roz90]} and \cite{[Pe-Za07]} for SDEs and SPDEs with martingale noises. In fact in \cite{[Alh-AMO-10]} we derived necessary conditions for optimality of stochastic systems similar to (\ref{intr:spde}), but the result there describes the maximum principle only in a local form and requires moreover the convexity of the control domain $U .$ In the present work we shall derive the maximum principle in its global form for our optimal control problem and, in particular, we shall not require the convexity of $U .$ Moreover, our results here generalize those in \cite{[Zhou93]} and \cite{[Be83]} and can be applied to the optimal control problem of partial observations with a given general nonlinear cost functional as done particularly in \cite[Section~6]{[Zhou93]}. The idea of reducing such a control problem to a control problem for a linear SPDE (Zakai's equation) was discussed also there. This is similar to (\ref{intr:spde}).

The main new features here are the driving noise is allowed to be an infinite dimensional martingale (as in Tudor \cite{[Tud89]} and Al-Hussein \cite{[Alh-Max10]}), the control domain $U$ need not be convex, and the control variable itself is allowed to enter in the martingale part of the equation as in the SPDE~(\ref{intr:spde}).

The present paper is organized as follows. In Section~\ref{sec2}
we introduce some definitions and notation that will be used throughout the paper.
In Section~\ref{sec3} our main stochastic control problem is introduced.
Section~\ref{sec4} is devoted to the adjoint equation of the SPDE~(\ref{intr:spde}) as well as the existence and uniqueness of its solution.
Finally, we state and establish the proof of our main result in Section~\ref{sec:main results}.

\section{Basic definitions and Notation}\label{sec2}
We assume that $(\Omega , \mathcal{F} , \{\mathcal{F}_t
\}_{t \geq 0 } , \mathbb{P} )$ is our complete filtered probability
space, such that $\{\mathcal{F}_t \}_{t \geq 0 } $ is a continuous filtration, in the sense that every square integrable $K$-valued
martingale with respect to $\{\mathcal{F}_t \, , \; 0 \leq t \leq T \}$ has a continuous version.
Let $\mathcal{P}$ denote the predictable $\sigma$\,-\,algebra of subsets of $\Omega\times
[0 , T] .$ A $K$\,-\,valued process is said to be predictable if it is
$\mathcal{P}/\mathcal{B}(K)$ measurable.
Let $\mathcal{M}^{2 , c}_{[0 , T ]} (K) $ be the space of all square
integrable continuous martingales in $K .$ We say that two
elements $M$ and $N$ of $\mathcal{M}^{2 , c}_{[0 , T]} (K) $ are {\it very strongly orthogonal} (VSO) if ${ \mathbb{E}\, [ M
(\tau )\otimes N (\tau ) ] =\, \mathbb{E}\, [ M (0)\otimes N (0) ] ,}$ for
all $[ 0 , T]$\,-\,valued stopping times $\tau .$

For $M \in \mathcal{M}^{2 , c}_{[0 , T ]} (K) $ let $<<M>>$ be its \emph{angle process} taking its values
in the space $L_1(K) ,$ where $L_1(K)$ is the space of nuclear operators on $K ,$ and satisfying ${ M\otimes M - <<M>> \, \in \mathcal{M}^{2 , c}_{[0 , T ]} (L_1 (K)) ,}$ and denote by $<M>$ the quadratic variation of $M .$ It is known (see \cite{[Me-P]}) that there exist a predictable process $\tilde{\mathcal{Q}}_{M} (s, \omega )$ in $L_1 (K)$ such that $<< M >>_t \; = \int_0^t \tilde{\mathcal{Q}}_{M} (s, \omega ) \, d <M>_s .$

For $(t , \omega )$ if $\tilde{\mathcal{Q}}(t , \omega )$ is any symmetric,
positive definite nuclear operator on $K ,$ we shall denote by
$L_{\tilde{\mathcal{Q}}(t , \omega )} (K)$ the set of all linear (not
necessarily bounded) operators $\Phi$ which map $\tilde{\mathcal{Q}}^{1/2}(t ,
\omega ) (K)$ into $K $ such that $\Phi \tilde{\mathcal{Q}}^{1/2}(t , \omega )
\in L_2 (K) ,$ the space of all Hilbert-Schmidt operators from
$K$ into itself. The inner product and norm in
$L_2 (K)$ will be denoted respectively by $\big{<} \cdot , \cdot \big{>}_2 $ and $|| \cdot
||_2 .$

We recall that the stochastic integral $\int_0^{\cdot} \Phi (s) d M(s) $ is defined
for mappings $\Phi$ such that for each $(t , \omega ) , \; \Phi (t ,
\omega ) \in L_{\tilde{\mathcal{Q}}_{_{M}} (t , \omega )} (K) ,$ for every $h \in K
$ the $K$\,-\,valued process $\Phi \tilde{\mathcal{Q}}^{1/2}_M
(h) $ is predictable, and
\( \mathbb{E}\; [ \, \int_0^T || ( \Phi \tilde{\mathcal{Q}}^{1/2}_M ) (t) ||_2^2
\; d <M>_t \, ] < \infty . \)

The space of such integrands is a Hilbert space with respect to the
scalar product $( \Phi_1 , \Phi_2 ) \mapsto \mathbb{E}\; [ \,
\int_0^T \big{<} \Phi_1 \tilde{\mathcal{Q}}^{1/2}_M \; ,
\Phi_2 \tilde{\mathcal{Q}}^{1/2}_M \big{>} \; d <M>_t \, ] .$ Simple processes in $L (K)$
are examples of such integrands. Hence the closure of the set of simple processes in this Hilbert space is
itself a Hilbert subspace. We denote it as in \cite{[Me-P]} by $\Lambda^2 ( K ; \mathcal{P} , M ) .$
More details and proofs can be found in \cite{[Me82]} or \cite{[Me-P]}.

In this paper we shall assume that there exists a measurable mapping ${ \mathcal{Q} (\cdot) : [ 0 , T ] \times \Omega \rightarrow L_1 (K)}$
such that $\mathcal{Q} (t)$ is symmetric,
positive definite, ${ << M >>_t \; = \int_0^t \mathcal{Q} (s) \, ds ,}$ and
$\mathcal{Q} (t) \leq \mathcal{Q} $ for some positive definite
nuclear operator $\mathcal{Q}$ on $K .$ Thus
$\tilde{\mathcal{Q}}_M (t) = \frac{\mathcal{Q} (t)}{q (t) }$ and
$<M>_t \, = \int_0^t q(s) ds ,$ with $q(t) = {\rm tr } \, (
\mathcal{Q} (t) ) .$ Thus, if $\Phi \in \Lambda^2 ( K ;
\mathcal{P} , M ) ,$
\begin{eqnarray*}
\mathbb{E} \, \Big[ \, | \int_0^T \Phi (s) dM(s) |^2 \Big] &=&
\mathbb{E} \, \Big[ \, \int_0^T || \Phi (s) \mathcal{Q}^{1/2} (s)
||_2^2 \, ds \Big] .
\end{eqnarray*}
This equality will be used frequently in the proofs given in Section~\ref{sec:main results}.
The process $\mathcal{Q} (\cdot)$ will play an essential role in deriving the adjoint equation of the SPDE~(\ref{intr:spde}), as appearing in the equation
(\ref{eq:defn of H}) in Section~\ref{sec4}; see in particular the discussion following equation~(\ref{eq:adjoint-bspde}).

\section{Statement of the control problem}\label{sec3}
Let us consider the following space:
\[ L^2_{\mathcal{F}} ( 0 , T ;  E  ) := \{ \psi : [ 0 , T]\times \Omega \rightarrow E,
\; \text{predictable and} \; \mathbb{E} \, [ \int_0^T | \psi (t) |^2 d
t \, ] < \infty \, \} ,\] where $E$ is a separable Hilbert
space. Suppose that $\mathcal{O}$ is a separable Hilbert space with an inner product $\big{<} \cdot , \cdot \big{>}_{\mathcal{O}}$, and $U$ is a nonempty subset of $\mathcal{O} .$ Denote by $$\mathcal{U}_{ad} = \{ u (\cdot ) : [0 , T]\times \Omega \rightarrow \mathcal{O}  \,\text{s.t.} \, u (\cdot ) \in  L^2_{\mathcal{F}} ( 0 , T ; \mathcal{O} )  \} .$$
This set is called the \emph{set of admissible controls} and its elements are called \emph{admissible controls}.

\bigskip

Now let us recall our SPDE:
\begin{eqnarray}\label{eq:see1}
 \left\{ \begin{array}{ll}
              d x(t) = ( A (t) x(t) + a (t , u (t) ) x(t) + b(t , u(t)) ) d t
              \\ \hspace{3cm} + \, [ \big< \sigma (t , u(t)) , x(t) \big>_{K} + \, g (t , u(t)) ] d M(t) , \\
             \; x(0) = x_0 \in K ,
         \end{array}
 \right.
\end{eqnarray}
and impose on it the following assumptions:

\bigskip

\noindent (i)\;  $A ( t , \omega )$ is a linear operator on $K ,$
$\mathcal{P}$\,-\,measurable, belongs to $L ( V ; V' ) $ uniformly
in $( t , \omega )$ and satisfies the following two conditions.
\begin{itemize}
\item{(1)} \, $A ( t , \omega )$ satisfies the coercivity
condition:
\[ 2 \; \big{<} A ( t , \omega ) \, y \, , y \big{>} + \; \alpha\; | y |^2_{_{V}}
\leq \lambda \; | y |^2  \; \; \;  a.e. \; t \in [ 0 , T ] \, , \; \
a.s. \; \; \forall \; y \in V , \] for some $\alpha , \lambda
> 0 .$ \\
\item{(2)} \, $\exists \; k_1 \geq 0 $ such that for all $( t , \omega ) $
\[ | \, A ( t , \omega ) \, y \, |_{_{V'}} \leq k_1 \, | y |_{_{V}} \; \;  \forall \; y \in V . \]
\end{itemize}
\noindent (ii) \; $a : \Omega \times [0, T] \times \mathcal{O} \rightarrow \mathbb{R}$ ,
$b : \Omega \times [0, T] \times \mathcal{O} \rightarrow K ,$
$\sigma : \Omega \times [0, T] \times \mathcal{O} \rightarrow K $ and
$g : \Omega \times [0, T] \times \mathcal{O} \rightarrow L_{\mathcal{Q}} (K) $
are predictable and bounded given mappings.

\bigskip

\begin{defn}\label{defn:soln of SEE}
We say that $x = x^{u (\cdot ) } \in L^2_{\mathcal{F}} ( 0 , T ;  V )$ is a \emph{solution} of (\ref{eq:see1}) if $\forall \; \eta \in V$ (or any dense subset) and for almost all $(t , \omega ) \in [0 , T] \times \Omega$
\begin{eqnarray*}
              \big< x(t) , \eta \big>_K &=& \big< x_0 , \eta \big>_K + \int_0^t \big< A (s) x(s) + a (s , u (s) ) x(s) + b(s , u(s)) \, , \eta \big>_{V} \, ds
              \\ && + \int_0^t \big< \eta \, , [ \big< \sigma (s , u(s)) , x(s) \big>_{K} + \, g (s , u(s)) ] d M(s) \big>_K . \\
\end{eqnarray*}
\end{defn}

\bigskip

Given a bounded measurable mapping $\ell : [ 0 , T] \times \mathcal{O} \rightarrow K$ and a fixed element $G$ of $K,$ we define the \emph{cost functional} by:
\begin{equation}\label{cost functional}
J(u (\cdot ) ) : = \mathbb{E} \; \big{[} \; \int_0^T \big< \ell (t , u (t) ) , x^{u
(\cdot ) } (t) \big>_{K} \, dt + \big< G , x^{u (\cdot ) } (T) \big>_{K} \; \big{]} , \; u (\cdot ) \in \mathcal{U}_{ad} .
\end{equation}

It is easy to realize that under assumptions (i) and (ii) there exists a unique solution to (\ref{eq:see1}) in $L^2_{\mathcal{F}} ( 0 , T ; K ) .$ This fact can be found in \cite[Theorem~4.1, P. 105]{[Gre-Tu95]}, \cite[Theorem~2.10]{[Gyo83]} or \cite[Theorem~3.2]{[Alh-Max10]}, and also can be gleaned from \cite{[Roz90]}. It\^{o}'s formula for such SPDEs can be found in \cite[Theorems~1, 2]{[Gyo-Kr81]}.

\smallskip

Our control problem is to minimize (\ref{cost functional}) over $\mathcal{U}_{ad} .$
Any $u^{*} ( \cdot ) \in \mathcal{U}_{ad} $ satisfying
\begin{equation}\label{value-function}
J(u^{*} ( \cdot ) ) = \inf \{ J( u (\cdot ) ): \; u (\cdot ) \in \,
\mathcal{U}_{ad} \}
\end{equation}
is called an \emph{optimal control}. The corresponding solution $x^{u^{*} (\cdot ) }$ of (\ref{eq:see1}), which we denote briefly by $x^{*}$ and $( x^{*} \, , u^{*} (\cdot ) )$ are called respectively an \emph{optimal solution} and an \emph{optimal pair} of the stochastic optimal control problem (\ref{eq:see1})-(\ref{value-function}).

\bigskip

The existence problem of optimal control can be developed from the works of \cite{[Ahmed86-existence]}, \cite{[Ahmed91-relaexed]} and \cite{[Tud89]}. However, a special case can be found in \cite{[Alh-Max10]}.

\section{Adjoint equation}\label{sec4}
Recall the SPDE~(\ref{eq:see1}) and the mappings in (\ref{cost functional}), and define the \emph{Hamiltonian}
${ H:[ 0 , T ] \times \Omega \times K \times \mathcal{O} \times K \times L_2 (K) \rightarrow
\mathbb{R} }$ for $( t , \omega , x , v , y , z ) \in [ 0 , T ] \times \Omega \times K \times \mathcal{O}
\times K \times L_2 (K) $ by
\begin{eqnarray}\label{eq:defn of H}
 H ( t , \omega , x , v , y , z ) &:=& - \big< \ell (t , v ) \, , x \big>_{V} - a (t , \omega , v) \big<
 x \, , y \big>_K \nonumber \\
&& - \, \big< b (t , \omega , v) \, , y \big>_K - \big< \tilde{\sigma} (t , \omega , x , v) \mathcal{Q}^{1/2} (t , \omega ) \, , z
\big>_2 \, ,
\end{eqnarray}
where $\tilde{\sigma} : [ 0 , T ] \times \Omega \times K  \times \mathcal{O} \rightarrow L_{\mathcal{Q}} (K) $ is defined by
\begin{eqnarray*}
\tilde{\sigma} (t , \omega , x , v) = \big< \sigma (t , \omega , v) \, , x \big>_K \, \Phi (x) + g (t, \omega , v)
\end{eqnarray*}
with $\Phi$ being the constant mapping $\Phi : K \rightarrow L_{\mathcal{Q}} (K) , x \mapsto \Phi (x) = \text{id}_K .$
Then
\begin{eqnarray*}
&&  \hspace{-0.5cm} \big< \tilde{\sigma} (t , \omega , x , v) \mathcal{Q}^{1/2} (t , \omega ) \, , z
\big>_2 =
\big< \, \big<  \sigma (t , \omega , v) , x \big>_K \; ( \Phi (x) + g(t , \omega , v) ) \mathcal{Q}^{1/2} (t , \omega ) \, , z \big>_2 \nonumber
\\ && \hspace{2cm} = \, \big< \, \big< \mathcal{Q}^{1/2} (t , \omega )  \, , z \big>_2 \, \sigma (t , \omega , v) , x \big>_K  + \big< g(t , \omega , v) \mathcal{Q}^{1/2} (t , \omega ) \, , z \big>_2
 \nonumber
\\ && \hspace{2cm} = \, \big< B (t , \omega , v ) z , x \big>_K + \big< g(t , \omega , v) \mathcal{Q}^{1/2} (t , \omega ) \, , z \big>_2 \, ,
\end{eqnarray*}
where $B : [0 , T] \times \Omega \times \mathcal{O} \rightarrow L ( L_2 (K) , K) $ is defined such that
\begin{equation}\label{eq:formula of B}
B(t , \omega , v) z = \big< \mathcal{Q}^{1/2} (t , \omega ) \, , z \big>_2 \; \sigma (t , \omega , v) .
\end{equation}
Moreover,
\begin{equation}\label{eq:derivative of H}
\nabla_{x} H ( t , \omega , x , v , y , z ) = - \ell (t , v ) - a (t , \omega , v) y - B(t , \omega , v) z .
\end{equation}

The adjoint equation of (\ref{eq:see1}) is the following BSPDE:
\begin{eqnarray}\label{eq:adjoint-bspde}
 \left\{ \begin{array}{ll}
             d y^{u (\cdot )}(t) = - \Big[  A^* (t) y(t) - \nabla_{x} H
             ( t , x^{u (\cdot )} (t), u (t), y^{u (\cdot )} (t) , z^{u (\cdot )}
              (t) \mathcal{Q}^{1/2} (t) ) \Big] \,
              dt  \\  \hspace{2.1in} + z^{u (\cdot )} (t)\, d M(t) + d N^{u (\cdot )} (t)
             , \; \; \; 0 \leq t < T , \\
             \;  y^{u (\cdot )} (T) = G ,
         \end{array}
 \right.
\end{eqnarray}
where $A^*(t)$ is the adjoint operator of $A(t) .$

\smallskip

It is important to realize that the presence of the process $\mathcal{Q}^{1/2} (\cdot )$ in the equation~(\ref{eq:adjoint-bspde}) is
crucial in order for the mapping $\nabla_{x} H$ to be defined on the space $L_2 (K) ,$ since the process $z^{u (\cdot )}$ need not be bounded as it is
discussed in Section~\ref{sec2}. This has to be taken always into account when dealing with BSPDEs and even BSDEs in infinite dimensions; cf. also \cite{[Alh-BSDE-2011]}.

\bigskip

The following theorem gives the solution to this BSPDE~(\ref{eq:adjoint-bspde}) in the sense that there exists a triple
$( y^{u (\cdot )} , z^{u (\cdot )} , N^{u (\cdot )} )$ in $L^2_{\mathcal{F}} ( 0 , T ; K )\times \Lambda^2 ( K ; \mathcal{P} ,
M ) \times \mathcal{M}^{2 , c}_{[ 0 , T] } (K) $ such that the
following equality holds $a.s. $ for all $t \in [ 0 , T ] , \; N (0)
= 0$ and $N$ is VSO to $M$:
\begin{eqnarray*}
y^{u (\cdot )} (t) &=& \xi + \int_t^T \nabla_{x} H
             ( s , x^{u (\cdot )} (s), u (s), y^{u (\cdot )} (s) , z^{u (\cdot )}
              (s) \mathcal{Q}^{1/2} (s) ) \,
              ds  \nonumber \\
&& - \int_t^T z^{u (\cdot )} (s ) \; d M (s) - \int_t^T d N^{u (\cdot )} (s) .
\end{eqnarray*}
\begin{thm}\label{th:solution of adjointeqn}
Assume that (i)--(ii) hold. Then there exists a unique solution
$( y^{u (\cdot )} , z^{u (\cdot )} , N^{u (\cdot )} )$ of the
BSDE~(\ref{eq:adjoint-bspde}) in $L^2_{\mathcal{F}} ( 0 , T ; K )\times
\Lambda^2 ( K ; \mathcal{P} , M ) \times \mathcal{M}^{2 , c}_{[ 0 ,
T] } (K) .$
\end{thm}
The proof of this theorem can be found in \cite{[Alh-Stoc09]}.

\smallskip

We shall denote briefly the solution of (\ref{eq:adjoint-bspde}) corresponding to the optimal control $u^{*} ( \cdot )$ by $(y^{*} , z^{*} , N^{*}) .$

\section{Main results}\label{sec:main results}
In this section we shall derive and prove our main result on the maximum principle for optimal control of the SPDE~(\ref{eq:see1}) associated with cost functional (\ref{cost functional}) and value function (\ref{value-function}) by using the results of the previous section on the adjoint equation (BSPDE). Before doing so, let us mention that the relationship between BSPDEs and maximum principle for some SPDEs is developed in several works, among them for instance are \cite{[Pe-93]} and \cite{[Y-Z]} and the references of Zhou cited therein. Other discussions in this respect can be found in \cite{[Tang-Li]} and \cite{[Zhou93]} as well. Bensoussan in \cite[Chapter~8]{[Be-book]} presents a stochastic maximum principle approach to the problem of stochastic control with partial information treating a general infinite dimensional setting and the adjoint equation is derived also there. Another work on the maximum principle that is connected to BSDEs can be found also in \cite{[Bah-Mez05]}. For an expanded discussion on the history of maximum principle we refer the reader to \cite[P. 153--156]{[Y-Z]}. And finally, one can find also useful information in Bensoussan's lecture notes \cite{[Be82]}, \cite{[Be82]} and Li \& Yong \cite{[Li-Yong95]} in addition to the references therein.

\bigskip

\pagebreak

Our main theorem is the following.
\begin{thm}\label{thm:main thm}
Suppose (i)--(ii). If $( x^{*} , u^{*} (\cdot ) )$ is an optimal
pair for the problem (\ref{eq:see1})-(\ref{value-function}), then there exists a unique solution $( y^{*},
z^{*} , N^{*} )$ to the corresponding BSEE~(\ref{eq:adjoint-bspde}) such that the following inequality holds:
\begin{eqnarray}\label{ineq1:main}
&& \hspace{-1cm} H ( t , x^{*} (t ) , u , y^{*} (t ) , z^{*} (t ) \mathcal{Q}^{1/2} (t ) )
 \nonumber \\ && \hspace{1.5cm} \leq \, H ( t , x^{*} (t ) , u^{*} (t ) , y^{*} (t) , z^{*} (t ) \mathcal{Q}^{1/2} (t ) )  \\
&& \hspace{2.5in}  \text{a.e.} \; t \in [0 , T],\; \text{a.s.} \; \forall \; u \in U .  \nonumber
\end{eqnarray}
\end{thm}

To start proving the theorem we need to develop some necessary estimates using the so-called spike variation method. For this we let $(x^{*} , u^{*} (\cdot ) )$ be the given optimal pair. Let $0 \leq t_0 < T$ be fixed such that $\mathbb{E} \; [ | x(t_0) |_K^2 ] < \infty $
 and $0 \leq \varepsilon < T - t_0 .$ Let $u$ be a random variable taking its values in $U ,$  $\mathcal{F}_{t_0}$\,-\,measurable and $\displaystyle{\sup_{\omega \in \Omega }} \, | u (\omega ) | < \infty .$ Consider the following spike variation of the control $u^{*} (\cdot )$:
\begin{eqnarray}\label{eq:defn of control var}
 u_{\varepsilon } (t) = \left\{ \begin{array}{ll}
u^{*} (t) \;\; & \text{ if} \; \; t \in [0 , T] \backslash [t_0 , t_0 + \varepsilon ] \nonumber \\
u \;\; & \text{ if} \; \; t \in [t_0 , t_0 + \varepsilon ] .
\end{array}
 \right.
\end{eqnarray}

We can consider the $x^{u_{\varepsilon } (\cdot ) }$ as the solution of the SPDE~(\ref{eq:see1}) corresponding to $ u_{\varepsilon } (\cdot ) .$ We shall denote it briefly by $x_{\varepsilon } .$ Note that $x_{\varepsilon } (t) = x^{*} (t) $ for all $0 \leq t \leq t_0 .$

\bigskip

We shall divide the proof into several lemmas as follows.

\begin{lem}\label{lem:lemma1}
Suppose (i)--(ii). Then
\begin{equation}\label{ineq:ineq1-lem1}
\sup_{t_0 \leq t \leq t_0 + \varepsilon } \mathbb{E} \, [ | x_{\varepsilon} (t) |_K^2 \, ] \, \leq C_1 \big( \, \mathbb{E} \, [ | x^* (t_0) |^2_K + C_2 \, \varepsilon \, \big)
\end{equation}
for some positive constants $C_1$ and $C_2 .$
\end{lem}
\begin{proof}
Observe first from (\ref{eq:see1}) and (\ref{eq:defn of control var}) that, for $t_0 \leq t \leq t_0 + \varepsilon ,$
\begin{eqnarray}\label{eq:eq2-lem1}
               x_{\varepsilon} (t) = x^{*} (t_0) + \int_{t_0}^t \big( A (s) x_{\varepsilon}(s) + a (s , u ) x_{\varepsilon}(s) + b(s , u) \big) d s
              \nonumber \\ \hspace{3cm} + \, \int_{t_0}^t \big[ \big< \sigma (s , u) , x_{\varepsilon}(s) \big>_{K} + \, g (s , u) \big] d M(s) .
\end{eqnarray}
Therefore, by It\^{o}'s formula, assumption (i), Cauchy-Schwartz inequality and assumption (ii) we get
\begin{eqnarray}\label{ineq:ineq3-lem1}
& & \hspace{-1cm} \mathbb{E}\, [ \, |\, x_{\varepsilon} (t) \, |^2_K \, ] + \alpha\; \mathbb{E}\,  [ \, \int_{t_0}^t | \, x_{\varepsilon} (s) \, |^2_{V} \, ds ] \leq \mathbb{E}\, [ \, |\, x^{*} (t_0) \, |^2_K \, ] + \lambda \, \mathbb{E}\,  [ \, \int_{t_0}^t | \, x_{\varepsilon} (s) \, |_K^2 \, ds \, ]
\nonumber \\
&& \hspace{-0.75cm} + \; 2 \, \mathbb{E}\,  [ \, \int_{t_0}^t \big< a (s , u) \, x_{\varepsilon} (s) \, , x_{\varepsilon} (s) \big>_K ds \, ]
+ 2 \, \mathbb{E}\,  [ \, \int_{t_0}^t \big<  x_{\varepsilon} (s) \, , b(s, u) \big>_K ds \, ]
\nonumber \\
&& \hspace{-0.75cm}
+ \; 2 \, \mathbb{E}\,  [ \, \int_{t_0}^t || \, \big< \sigma (s , u) \, , x_{\varepsilon } (s) \big>_{K} \, \text{id}_K \, \mathcal{Q}^{1/2} (s) \, ||^2_2 ds \, ]
+ 2 \, \mathbb{E}\,  [ \, \int_{t_0}^t || \, g(s , u) \mathcal{Q}^{1/2} (s) \, ||^2_2 ds \, ]
\nonumber \\ &\leq&
\mathbb{E}\, [ \, |\, x^{*} (t_0) \, |^2_K \, ] + \lambda \, \mathbb{E}\, [ \, \int_{t_0}^t | \, x_{\varepsilon} (s) \, |_K^2 \, ds \, ] +  2\, k_1 \, \mathbb{E}\,  [ \, \int_{t_0}^t | \, x_{\varepsilon} (s) \, |_K^2 \, ds \, ] \nonumber \\
&& \hspace{2in} + \; k_2^2 \, \mathbb{E}\,  [ \, \int_{t_0}^t | \, x_{\varepsilon} (s) \, |_K^2 \, ds \, ] + (t-t_0) \nonumber \\
&& \hspace{4cm} + \; 2 k_3^2 \, || \mathcal{Q}^{1/2} ||^2_2 \; \mathbb{E}\,  [ \, \int_{t_0}^t | \, x_{\varepsilon} (s) \, |_K^2 \, ds \, ]  +  2 k_4^2 \varepsilon \, || \mathcal{Q}^{1/2} ||^2_2
\nonumber \\ &=& \big( \lambda + 2 k_1 + k_2^2 + 2 k_3^2 \, || \mathcal{Q}^{1/2} ||^2_2 \; (1+ \varepsilon) \big) \;  \int_{t_0}^t \mathbb{E}\, [ | \, x_{\varepsilon} (s) |_K^2 \, ] \, ds
\nonumber \\
&& \hspace{2in} + \; ( 1 + 2 k_4 k_4^2 \, || \mathcal{Q}^{1/2} ||^2_2 ) \, \varepsilon + \mathbb{E}\, [ \, |\, x^{*} (t_0) \, |^2_K \, ] .
\end{eqnarray}
In the last part of this inequality we have used the boundedness in assumption (ii) of the mappings $a , b , \sigma , g$ respectively to get the constants $k_1 - k_4 .$

Thus, in particular, by applying Gronwall's inequality to (\ref{ineq:ineq3-lem1}) we obtain (\ref{ineq:ineq1-lem1}) with
$$ C_1 = e^{\varepsilon \, \big( \lambda + 2 k_1 + k_2^2 + 2 k_3^2 \, || \mathcal{Q}^{1/2} ||^2_2 (1+ \varepsilon) \big)) } $$
and
$$ C_2 = 1 + 2 k_4^2 \, || \mathcal{Q}^{1/2} ||^2_2 .$$
This completes the proof.
\end{proof}

\begin{lem}\label{lem:lemma2}
Suppose (i)--(ii). Then
\begin{equation}\label{ineq:ineq1-lem2}
\sup_{t_0 + \varepsilon \leq t \leq T } \mathbb{E} \, [ | x_{\varepsilon} (t) |_K^2 \, ] \, \leq C_3 \big( \, \mathbb{E} \, [ | x^* (t_0) |^2_K ] + C_4 \, \varepsilon +1 \, \big)
\end{equation}
for some positive constants $C_3$ and $C_4 .$
\end{lem}
\begin{proof}
For $t_0 + \varepsilon \leq t \leq T,$ it follows that
\begin{eqnarray}\label{eq:eq1-lem2}
               x_{\varepsilon} (t) = x^{*} (t_0 + \varepsilon ) + \int_{t_0 + \varepsilon}^t \big( A (s) x_{\varepsilon}(s) + a (s , u^{*} (s) ) x_{\varepsilon}(s) + b(s , u^{*} (s) ) \big) d s
              \nonumber \\ \hspace{3cm} + \, \int_{t_0 + \varepsilon}^t \big[ \big< \sigma (s , u^{*} (s) ) , x_{\varepsilon}(s) \big>_{K} + \, g (s , u^{*} (s) ) \big] d M(s) .
\end{eqnarray}

Thus mimicking the proof of Lemma~\ref{lem:lemma1} and then applying inequality (\ref{ineq:ineq1-lem1}) easily yields (\ref{ineq:ineq1-lem2}).
\end{proof}

\begin{lem}\label{lem:lemma3}
Suppose (i)--(ii). Let $\xi_{\varepsilon} (t) = x_{\varepsilon} (t) - x^{*} (t), $ for $t \in [0 , T] .$ Then
\begin{equation}\label{eq:eq1-lem3}
\sup_{t_0 + \varepsilon \leq t \leq T } \mathbb{E} \, [ |\xi_{\varepsilon} (t) |_K^2 \, ] = O ( \varepsilon ) .
\end{equation}
\end{lem}
\begin{proof}
It is easy to get for $t \in [t_0 + \varepsilon , T] ,$
\begin{eqnarray}\label{eq:eq2-lem3}
\xi_{\varepsilon } (t) &=& \xi_{\varepsilon } (t_0 + \varepsilon ) + \int_{t_0 + \varepsilon}^t \big( A(s) \xi_{\varepsilon } (s) + a (s , u^* (s)) \, \xi_{\varepsilon} (s) \big) ds \nonumber \\ && \hspace{3cm} + \int_{t_0 + \varepsilon}^t \big< \sigma (s , u^* (s)) \, , \xi_{\varepsilon } (s) \big>_{K} \, dM(s) .
\end{eqnarray}
Hence, as done in the proof of Lemma~\ref{lem:lemma1}, we get
\begin{equation}\label{eq:eq3-lem3}
\sup_{t_0 + \varepsilon \leq t \leq T } \mathbb{E} \, [ |\xi_{\varepsilon} (t) |_K^2 \, ] \leq C_5  \, \mathbb{E} \, [ |\xi_{\varepsilon} (t_0 + \varepsilon) |_K^2 \, ] .
\end{equation}

On the other hand, for $t_0 \leq t \leq t_0 + \varepsilon$ we have $\xi_{\varepsilon } (t_0 ) = 0 $ and
\begin{eqnarray}\label{eq:eq4-lem3}
\xi_{\varepsilon } (t) &=& \int_{t_0}^t \Big[ A(s) \xi_{\varepsilon } (s) + \big( a(s , u) - a (s , u^* (s)) \big) \, x_{\varepsilon} (s) \nonumber \\ && \hspace{1.5cm}
+ \; \big( b (s , u) - b (s , u^* (s)) \big) + a (s , u^* (s)) \, \xi_{\varepsilon} (s) \Big] ds \nonumber \\
&& + \; \int_{t_0}^t \Big[ \big< \sigma (s , u) - \sigma (s , u^* (s)) \, , x_{\varepsilon } (s) \big>_{K} \nonumber \\
&& +  \; \big( g (s , u) - g (s , u^* (s)) \big) + \big< \sigma (s , u^* (s)) \, , \xi_{\varepsilon } (s) \big>_{K} \Big] dM(s) .
\end{eqnarray}
Hence by It\^{o}'s formula, assumption (i), Cauchy-Schwartz inequality and assumption (ii) it follows that
\begin{eqnarray}\label{ineq:ineq5-lem3}
& & \hspace{-1cm} \mathbb{E}\, [ \, |\, \xi_{\varepsilon} (t) \, |^2_K \, ] + \alpha\; \mathbb{E}\,  [ \, \int_{t_0}^t | \, \xi_{\varepsilon} (s) \, |^2_{V} \, ds ] \nonumber \\ &\leq& \lambda \, \mathbb{E}\,  [ \, \int_{t_0}^t | \, \xi_{\varepsilon} (s) \, |_K^2 \, ds \, ] + \,
 2 \, \mathbb{E} \, \big[ \, \int_{t_0}^{t} \big< \xi_{\varepsilon} (s) \, , \big( a(s , u) - a (s , u^* (s)) \big) \, x_{\varepsilon} (s) \big>_{K} \, ds \, \big]  \nonumber \\
 && \hspace{1cm} + \, 2 \, \mathbb{E} \, \big[ \, \int_{t_0}^{t} \big< \xi_{\varepsilon} (s) \, , b (s , u) - b (s , u^* (s)) \big>_{K} \, ds \, \big]
\nonumber \\
&& \hspace{1cm} + \, 2 \, \mathbb{E} \, \big[ \, \int_{t_0}^{t} \big< \xi_{\varepsilon} (s) \, , a (s , u^* (s)) \, \xi_{\varepsilon} (s) \big>_{K} \, ds \, \big]
\nonumber \\
&& \hspace{1cm} + \, 3 \, \mathbb{E}\,  [ \, \int_{t_0}^t || \, \big< \sigma (s , u)  -  \sigma (s , u^* (s)) \, , x_{\varepsilon } (s) \big>_{K} \, \text{id}_K  \, \mathcal{Q}^{1/2} (s) \, ||^2_2 ds \, ] \nonumber \\
&& \hspace{1cm} + \, 3 \, \mathbb{E}\,  [ \, \int_{t_0}^t || \, \big< \sigma (s , u^* (s)) \, , \xi_{\varepsilon } (s) \big>_{K} \, \text{id}_K  \, \mathcal{Q}^{1/2} (s) \, ||^2_2 ds \, ]
\nonumber \\
&& \hspace{1cm}
+ 3 \, \mathbb{E}\,  [ \, \int_{t_0}^t || \, \big( g(s , u) - g(s , u^* (s)) \big) \, \mathcal{Q}^{1/2} (s) \, ||^2_2 ds \, ]
\nonumber \\
&\leq& ( \lambda + 4 k_1^2 + k_2^2 + 2 k_1 + 3 k_3^2 \cdot || \mathcal{Q}^{1/2} ||^2_2 ) \,  \int_{t_0}^t \mathbb{E}\, [ | \, \xi_{\varepsilon} (s) |_K^2 \, ] \, ds
\nonumber \\
&& + \, ( 6 k_3^2 \, || \mathcal{Q}^{1/2} ||^2_2 + 1 ) \,  \int_{t_0}^{t_0 + \varepsilon} \mathbb{E}\, [ | \, x_{\varepsilon} (s) |_K^2 \, ] \, ds + (1 + 12 k_4^2 ) \varepsilon
\nonumber \\
&\leq& ( \lambda + 4 k_1^2 + k_2^2 + 2 k_1 + 3 k_3^2 \, || \mathcal{Q}^{1/2} ||^2_2 ) \,  \int_{t_0}^t \mathbb{E}\, [ | \, \xi_{\varepsilon} (s) |_K^2 \, ] \, ds
\nonumber \\
\nonumber \\
&& + \, ( 6 k_3^2 \, || \mathcal{Q}^{1/2} ||^2_2 + 1 ) \, C_1 \cdot \big(  \mathbb{E} \, [ | \, x^{*} (t_0) |_K^2 \, ] + C_2 \, \varepsilon \big) \, \varepsilon + (1 + 12 k_4^2 ) \varepsilon .
\end{eqnarray}

Therefore Gronwall's inequality gives
\begin{equation}\label{ineq:ineq6-lem3}
\sup_{t_0 \leq t_0 + \varepsilon } \; \mathbb{E}\, [ \, |\, \xi_{\varepsilon} (t) \, |^2_K \, ]
\leq C_6 ( \varepsilon ) \cdot \varepsilon ,
\end{equation}
where
\begin{eqnarray*}
&& \hspace{-0.5cm} C_6 ( \varepsilon ) = e^{( \lambda + 4 k_1^2 + k_2^2 + 2 k_1 + 3 k_3^2 \, || \mathcal{Q}^{1/2} ||^2_2 ) \, \varepsilon } \cdot \Big[ \, ( 6 k_3^2 \, || \mathcal{Q}^{1/2} ||^2_2 + 1 ) \, C_1 \cdot \big( \mathbb{E} \, [ | \, x^{*} (t_0) |_K^2 \, ] \\ && \hspace{3.75in} + \; C_2 \, \varepsilon \big) + 1 + 12 k_4^2 \, \Big] .
\end{eqnarray*}
Now by applying (\ref{ineq:ineq6-lem3}) in (\ref{eq:eq3-lem3}) it yields eventually
\begin{equation}\label{eq:eq7-lem3}
\sup_{t_0 + \varepsilon \leq t \leq T } \mathbb{E} \, [ |\xi_{\varepsilon} (t) |_K^2 \, ] \leq C_5  \, C_6 (\varepsilon) \cdot \varepsilon .
\end{equation}
Thus (\ref{eq:eq1-lem3}) follows.
\end{proof}

\bigskip

In the following result we shall try to compute $\mathbb{E} \, [ \, \big< \, y^{*} (t_0 + \varepsilon ) , \xi (t_0 + \varepsilon ) \, \big>_K \, ] .$
\begin{lem}\label{lem:lemma4}
Suppose (i)--(ii). We have
\begin{eqnarray}\label{eq:Ito formula of Y and xi-1}
&& \hspace{-0.5cm} \mathbb{E} \, \Big{[} \, \big< \, y^{*} (t_0 + \varepsilon ) \, , \xi_{\varepsilon } (t_0 + \varepsilon ) \, \big>_K + \int_{t_0}^{t_0 + \varepsilon} \big< \ell (t, u^{*} (t) ) \, , \xi_{\varepsilon} (t) \big>_{K} dt \, \Big] \nonumber \\
&& \hspace{1cm} = \, \mathbb{E} \, \big{[} \, \int_{t_0}^{t_0 + \varepsilon} \big< y^* (t) \, , \big( a(t , u) - a (t , u^* (t)) \big) \, x_{\varepsilon} (t) \big>_{K} \, \big{]}  \nonumber \\
&& \hspace{0.5cm} + \, \mathbb{E} \, \big{[} \, \int_{t_0}^{t_0 + \varepsilon} \big< y^* (t) \, , b (t , u) - b (t , u^* (t)) \big>_{K} \, dt \, \big{]}
\nonumber \\
&& \hspace{0.5cm} + \, \mathbb{E} \, \big{[} \, \int_{t_0}^{t_0 + \varepsilon} \big< \sigma (t , u) - \sigma (t , u^* (t)) \, , x_{\varepsilon } (t) \big>_{K} \, \big<  \mathcal{Q}^{1/2} (t) \, , z^{*} (t)  \mathcal{Q}^{1/2} (t) \big>_2 \, dt \, \big{]}
\nonumber \\
&& \hspace{0.5cm} + \, \mathbb{E} \, \big{[} \, \int_{t_0}^{t_0 + \varepsilon} \big< \big( g (t , u) - g (t , u^* (t)) \big) \mathcal{Q}^{1/2} (t) \, , z^{*} (t)  \mathcal{Q}^{1/2} (t) \big>_2 \, dt \, \big{]}
\end{eqnarray}
and
\begin{eqnarray}\label{eq:Ito formula of Y and xi-2}
&& \hspace{-1.5cm} \mathbb{E} \, [ \, \big< \, y^{*} (t_0 + \varepsilon ) , \xi_{\varepsilon } (t_0 + \varepsilon ) \, \big>_K \, ] =  \mathbb{E} \, \big{[} \, \int_{t_0 + \varepsilon}^T \big< \ell (t, u^{*} (t) ) \, , \xi_{\varepsilon} (t) \big>_{K} dt \, \big{]} \nonumber \\
&& \hspace{3in} + \; \mathbb{E} \, \big{[} \, \big< G \, , \xi_{\varepsilon} (T) \big>_{K} \, \big{]} .
\end{eqnarray}
\end{lem}
\begin{proof}
Note that for $t_0 \leq t \leq t_0 + \varepsilon$ we have $\xi_{\varepsilon } (t_0 ) = 0$ and (\ref{eq:eq4-lem3}).
Therefore by using It\^o's formula to (\ref{eq:eq4-lem3}) together with (\ref{eq:adjoint-bspde}), (\ref{eq:derivative of H}) and (\ref{eq:formula of B}) we get easily (\ref{eq:Ito formula of Y and xi-1}).
The equality in (\ref{eq:Ito formula of Y and xi-2}) is proved similarly with the help of (\ref{eq:eq2-lem3}).
\end{proof}

\begin{lem}\label{lem:lemma5}
Suppose (i)--(ii). We have
\begin{eqnarray}\label{ineq1:main ineq}
&& \hspace{-0.75cm} 0 \leq \mathbb{E} \, \Big{[} \, \int_{t_0}^{t_0 + \varepsilon} \big< \ell (t, u ) - \ell (t, u^{*} (t) ) \, , x^{*} (t) \big>_{K} dt \, ]  \nonumber \\
&& \hspace{-0.35cm} + \, \mathbb{E} \, \big{[} \, \int_{t_0}^{t_0 + \varepsilon} \big< y^* (t) \, , \big( a (t , u) - a (t , u^* (t) ) \big) \, x^{*} (t) \big>_{K} \, \big{]}  \nonumber \\
&& \hspace{-0.35cm} + \, \mathbb{E} \, \big{[} \, \int_{t_0}^{t_0 + \varepsilon} \big< \sigma (t , u) - \sigma (t , u^* (t)) \, ,  \, x^{*} (t) \big>_{K} \, \big<  \mathcal{Q}^{1/2} (t) \, , z^{*} (t)  \mathcal{Q}^{1/2} (t) \big>_2 \, dt \, \big{]}
\nonumber \\
&& \hspace{-0.35cm} + \, \mathbb{E} \, \big{[} \, \int_{t_0}^{t_0 + \varepsilon} \big< y^* (t) \, , b (t , u) - b (t , u^* (t)) \big>_{K} \, dt \, \big{]}
\nonumber \\
&& \hspace{-0.35cm} + \, \mathbb{E} \, \big{[} \, \int_{t_0}^{t_0 + \varepsilon} \big< \big( g (t , u) - g (t , u^* (t) \big) \mathcal{Q}^{1/2} (t) \, , z^{*} (t)  \mathcal{Q}^{1/2} (t) \big>_2 \, dt \, \big{]} + o (\varepsilon ) .
\end{eqnarray}
\end{lem}
\begin{proof}
Since $u^{*} (\cdot )$ is optimal, we have
\begin{eqnarray*}
&& \hspace{-1.25cm} 0 \leq J ( u_{ \varepsilon } ( \cdot ) ) - J ( u^{*} ( \cdot ) ) \nonumber \\ &=&
 \mathbb{E} \; \big[ \; \int_0^T \big( \big< \ell (t , u_{ \varepsilon } (t) ) \, , x_{ \varepsilon } (t) \big>_{K} - \big< \ell (t , u^{*} (t) ) \, , x^{*} (t) \big>_{K} \big) \, dt \big] \nonumber
 \\ &&  \hspace{1.5cm} + \; \mathbb{E} \; \big[ \, \big< G \, , x_{ \varepsilon } (T) \big>_{K} - \big< G \, , x^{*} (T) \big>_{K} \, \big]
 \nonumber \\ &=&
\mathbb{E} \, \Big[ \, \int_{t_0}^{t_0 + \varepsilon } \big( \big< \ell (t , u) - \ell (t , u^{*} (t) ) \,, x_{ \varepsilon } (t) \big>_{K} + \big< \ell (t , u^{*} (t) ) \,, \xi_{ \varepsilon } (t) \big>_{K} \big) \, dt \Big] \nonumber
 \\ && \hspace{1.5cm} + \; \mathbb{E} \; \big[ \, \int_{t_0 + \varepsilon }^T \big< \ell (t , u^{*} (t) ) \,, \xi_{ \varepsilon } (t) \big>_{K} \, dt + \big< G , \xi_{ \varepsilon } (T) \big>_{K} \, \big] .
\end{eqnarray*}
Hence using Lemma~\ref{lem:lemma4} (\ref{eq:Ito formula of Y and xi-2}) in this inequality gives
\begin{eqnarray}\label{ineq:variational inequality-2}
&& 0 \leq
 \mathbb{E} \, \Big[ \, \int_{t_0}^{t_0 + \varepsilon } \big( \big< \ell (t , u) - \ell (t , u^{*} (t) ) \,, x_{ \varepsilon } (t) \big>_{K} dt + \big< \ell (t , u^{*} (t) ) \,, \xi_{ \varepsilon } (t) \big>_{K} \big) \, dt \Big] \nonumber
 \\ && \hspace{2.2in} +  \; \mathbb{E} \; \big[ \, \big< y^{*} (t_0 + \varepsilon ) ) \,, \xi_{\varepsilon } (t_0 + \varepsilon ) \big>_{K} \, \big] .
\end{eqnarray}

Again by Lemma~\ref{lem:lemma4} (\ref{eq:Ito formula of Y and xi-1}) inequality (\ref{ineq:variational inequality-2}) becomes
\begin{eqnarray}\label{ineq:variational inequality-3}
&& \hspace{-0.5cm} 0 \leq
 \mathbb{E} \; \big[ \, \int_{t_0}^{t_0 + \varepsilon } \big< \ell (t , u) - \ell (t , u^{*} (t) ) \,, x_{ \varepsilon } (t) \big>_{K} \, dt \, \big] \nonumber \\
&& \hspace{1cm} + \, \mathbb{E} \, \big{[} \, \int_{t_0}^{t_0 + \varepsilon} \big< y^* (t) \, , \big( a(t , u) - a (t , u^* (t)) \big) \, x_{\varepsilon} (t) \big>_{K}  \, dt \, \big{]}  \nonumber \\
&& \hspace{1cm} + \, \mathbb{E} \, \big{[} \, \int_{t_0}^{t_0 + \varepsilon} \big< y^* (t) \, , b (t , u) - b (t , u^* (t)) \big>_{K}  \, dt \, \big{]}
\nonumber \\
&& \hspace{1cm} + \, \mathbb{E} \, \big{[} \, \int_{t_0}^{t_0 + \varepsilon} \big< \sigma (t , u) - \sigma (t , u^* (t)) \, , x_{\varepsilon } (t) \big>_{K} \, \big<  \mathcal{Q}^{1/2} (t) \, , z^{*} (t)  \mathcal{Q}^{1/2} (t) \big>_2 \, dt \, \big{]}
\nonumber \\
&& \hspace{1cm} + \, \mathbb{E} \, \big{[} \, \int_{t_0}^{t_0 + \varepsilon} \big< \big( g (t , u) - g (t , u^* (t)) \big) \mathcal{Q}^{1/2} (t) \, , z^{*} (t)  \mathcal{Q}^{1/2} (t) \big>_2 \, dt \, \big{]} .
\end{eqnarray}

On the other hand, assumption (ii) and Lemma~\ref{lem:lemma3} imply
\begin{eqnarray}\label{ineq:variational inequality-4}
&& \hspace{-2cm} \frac{1}{\varepsilon } \; \mathbb{E} \, \big[ \, \int_{t_0}^{t_0 + \varepsilon} \big< y^* (t) \, , \big( a(t , u) - a (t , u^* (t)) \big) \, \xi_{\varepsilon} (t) \big>_{K}  \, dt \, \big]
\nonumber \\
&\leq&   C_7 \, (\frac{1}{\varepsilon }) \; \int_{t_0}^{t_0 + \varepsilon} \mathbb{E} \, \Big(  | y^* (t) |_K \cdot | \xi_{\varepsilon} (t) |_K \Big) \, dt
\nonumber \\
&\leq&   C_7 \, (\frac{1}{\varepsilon }) \; \int_{t_0}^{t_0 + \varepsilon} \Big( ( \frac{\varepsilon^{1/3}}{2}) \; \mathbb{E} \,  [ | y^* (t) |^2_K ]
 + (\frac{1}{2 \, \varepsilon^{1/3}}) \; \mathbb{E} \,  [ | \xi_{\varepsilon } (t) |^2_K ] \Big) \, dt
\nonumber \\
&\leq&   C_8 \, \Big( \, \varepsilon^{1/3} \, (\frac{1}{\varepsilon }) \; \int_{t_0}^{t_0 + \varepsilon} \; \mathbb{E} \,  [ | y^* (t) |^2_K ] dt
 + (\frac{1}{\varepsilon }) \, \varepsilon \, (\frac{1}{\varepsilon^{1/3}}) \, \varepsilon \, \Big)
\rightarrow 0 ,
\end{eqnarray}
as $\varepsilon \rightarrow 0 ,$ provided that $t_0$ is a Lebesgue point of the function $t \mapsto \mathbb{E} \, [ \, | y^{*} (t) \, |^2_K \, ] ,$  for some positive constants $C_7$ and $C_8 .$

Similarly,
\begin{eqnarray}\label{ineq:variational inequality-5}
&& \hspace{-2cm} \frac{1}{\varepsilon } \; \mathbb{E} \; \Big[ \, \int_{t_0}^{t_0 + \varepsilon } \Big( \big< \ell (t , u) - \ell (t , u^{*} (t) ) \,, \xi_{ \varepsilon } (t) \big>_{K} \nonumber \\
&&  \hspace{-1cm} + \, \big< \sigma (t , u) - \sigma (t , u^* (t)) \, , \xi_{\varepsilon } (t) \big>_{K} \, \big<  \mathcal{Q}^{1/2} (t) \, , z^{*} (t)  \mathcal{Q}^{1/2} (t) \big>_2 \Big) \, dt \, \Big] \rightarrow 0 ,
\end{eqnarray}
as $\varepsilon \rightarrow 0 ,$ provided that $t_0$ is a Lebesgue point of the function ${t \mapsto \mathbb{E} \, [ \, || z^{*} (t) \mathcal{Q}^{1/2} (t) \, ||^2_2 \, ] .}$

Therefore, by applying (\ref{ineq:variational inequality-4}) and (\ref{ineq:variational inequality-5}) in (\ref{ineq:variational inequality-3}) we obtain (\ref{ineq1:main ineq}).
\end{proof}

\bigskip

We are now ready to complete the proof of Theorem~\ref{thm:main thm}.

\noindent {\bf Proof of Theorem~\ref{thm:main thm} }
Divide (\ref{ineq1:main ineq}) in Lemma~\ref{lem:lemma5} by $\varepsilon$ and let $\varepsilon \rightarrow 0$ to get
\begin{eqnarray*}
&& \hspace{-0.55cm} \mathbb{E} \, \Big{[} \, \big< \ell (t_0, u) - \ell (t_0, u^{*} (t_0) ) \, , x^{*} (t_0) \big>_{K} + \big< y^* (t_0) \, , \big( a (t_0 , u) - a (t_0 , u^* (t_0) ) \big) \, x^{*} (t_0) \big>_{K} \, \big{]}  \nonumber \\
&& \hspace{1.45cm} + \, \mathbb{E} \, \big{[} \, \big< y^* (t_0) \, , b (t_0 , u) - b (t_0 , u^* (t_0)) \big>_{K} \, \big{]}
\nonumber \\
&& \hspace{1.45cm} + \, \mathbb{E} \, \big{[} \, \big< \sigma (t_0 , u) - \sigma (t_0 , u^* (t_0)) \, ,  \, x^{*} (t_0) \big>_{K} \, \big<  \mathcal{Q}^{1/2} (t_0) \, , z^{*} (t_0)  \mathcal{Q}^{1/2} (t_0) \big>_2 \, \big{]}
\nonumber \\
&&
\hspace{1.45cm} + \, \mathbb{E} \, \big{[} \, \big< \big( g (t_0 , u) - g (t_0 , u^* (t_0) \big) \mathcal{Q}^{1/2} (t_0) \, , z^{*} (t_0)  \mathcal{Q}^{1/2} (t_0) \big>_2 \, \big{]}
\geq 0 .
\end{eqnarray*}
Consequently,
\begin{eqnarray*}
&& \hspace{-1cm} H ( t_0 , x^{*} (t_0 ) , u , y^{*} (t_0 ) , z^{*} (t_0 ) \mathcal{Q}^{1/2} (t_0 ) )
 \nonumber \\ && \hspace{1.5cm} \leq \, H ( t_0 , x^{*} (t_0 ) , u^{*} (t_0 ) , y^{*} (t_0) , z^{*} (t_0 ) \mathcal{Q}^{1/2} (t_0 ) )  .
\end{eqnarray*}
Hence (\ref{ineq1:main}) holds by a standard argument as for example in \cite[Chapet~3]{[Y-Z]}, and the proof of Theorem~\ref{thm:main thm} is then complete.

\bigskip

{\bf Acknowledgement.}
This author would like to express his great thanks to Professor David Elworthy for reading the first draft of the paper and providing useful comments. Many thanks go to the Mathematics Institute, Warwick University, where this work was done,
for hospitality during the summer of 2011.

\fussy


\begin{thebibliography}{99}
\bibitem{[Ahmed86-existence]} N. U. Ahmed, Existence of optimal controls for a class of systems governed
by differential inclusions on a Banach space, J. Optim. Theory Appl., 50, No. 2 (1986), 213--237.
\bibitem{[Ahmed91-relaexed]} N. U. Ahmed, Relaxed controls for stochastic boundary value problems in infinite dimension.
Optimal control of partial differential equations (Irsee, 1990), 1--10, Lecture Notes in Control and Inform. Sci., 149, Springer, Berlin, 1991.
\bibitem{[Alh-Stoc09]} A. Al-Hussein, Backward stochastic partial
differential equations driven by infinite dimensional martingales
and applications, Stochastics, 81, 6 (2009) 601--626.
\bibitem{[Alh-Max10]} A. Al-Hussein, Maximum principle for controlled stochastic evolution
equations, Int. Journal of Math. Analysis, Vol. 4, No. 30 (2010), 1447--1464.
\bibitem{[Alh-AMO-10]} A. Al-Hussein, Necessary conditions for optimal control of stochastic evolution equations in
Hilbert spaces, Appl. Math. Optim. 63 (2011), No. 3, 385--400.
\bibitem{[Alh-BSDE-2011]} A. Al-Hussein, Pontryagin's maximum principle for optimal control of infinite dimensional SDEs, Preprint 2011.
\bibitem{[Bah-Mez05]} S. Bahlali and B. Mezerdi, A general
stochastic maximum principle for singular control problems,
Electron. J. Probab., 10, No. 30 (2005), 988--1004.
\bibitem{[Bal-Pard05]} V. Bally, V., \'{E}. Pardoux L. and Stoica, Backward
stochastic differential equations associated to a symmetric Markov
process, Potential Anal., 22, No. 1 (2005), 17--60.
\bibitem{[Be82]} A. Bensoussan, Lectures on stochastic control. Nonlinear filtering
and stochastic control (Cortona, 1981), 1--62, Lecture Notes in Math., 972, Springer, Berlin-New York, 1982.
\bibitem{[Be83]}  A. Bensoussan, Maximum principle and dynamic programming approaches of the optimal control of partially observed diffusions,
Stochastics 9 (1983), No. 3, 169--222.
\bibitem{[Be-book]} A. Bensoussan, Stochastic control of partially observable systems,
Cambridge University Press, Cambridge, 1992.
\bibitem{[Fuh-Tess]} M. Fuhrman and G. Tessitore, Nonlinear Kolmogorov
equations in infinite dimensional spaces: the backward stochastic
differential equations approach and applications to optimal control,
Ann. Probab., 30, No. 3 (2002), 1397--1465.
\bibitem{[Gre-Tu95]} W. Grecksch and C. Tudor, Stochastic
evolution equations. A Hilbert space approach. Mathematical Research
85, Akademie-Verlag, Berlin, 1995.
\bibitem{[Gyo-Kr81]} I. Gy\"{o}ngy and N. V. Krylov, On stochastics
equations with respect to semimartingales. II. It\^{o} formula in
Banach spaces, Stochastics, 6, No. 3-4 (1981/82), 153--173.
\bibitem{[Gyo83]} I. Gy\"{o}ngy, On stochastic equations with respect to semimartingales. III,
Stochastics, 7 (1982), No. 4, 231--254.
\bibitem{[Imk10]} P. Imkeller, A. Reveillac and A. Richter, Differentiability of quadratic BSDEs
generated by continuous martingales, arXiv:0907.0941 [math.PR], 2010.
\bibitem{[Kot84]} P. Kotelenez, A stopped Doob inequality for stochastic convolution integrals and stochastic
evolution equations, Stochastic Anal. Appl., 2, No. 3 (1984),
245--265.
\bibitem{[Kry-Roz07]}  N. V. Krylov and B. Rozovskii, Stochastic
evolution equations, in: Stochastic differential equations: theory
and applications, pp. 1--69, Interdiscip. Math. Sci. 2, World Sci.
Publ., NJ, Hackensack, 2007.
\bibitem{[Li-Yong95]} X. J. Li and J. M. Yong,
Optimal control theory for infinite-dimensional systems. Systems \& Control: Foundations \& Applications,
Birkhauser Boston, Inc., Boston, MA, 1995.
\bibitem{[Marie09]} M.-A. Morlais, Quadratic BSDEs driven by a continuous martingale
and applications to the utility maximization problem, Finance Stoch. 13 (2009),
No. 1, 121--150.
\bibitem{[Me82]} M. M\'etivier, Semimartingales. A course on stochastic
processes, de Gruyter Studies in Mathematics 2, Walter de Gruyter \&
Co., Berlin-New York, 1982.
\bibitem{[Me-P]} M. M\'etivier and J. Pellaumail, Stochastic
integration, Probability and Mathematical Statistics, Academic Press
[Harcourt Brace Jovanovich, Publishers], New York-London-Toronto,
1980.
\bibitem{[Oks-Z05]} B. {\O}ksendal, F. Proske and T. Zhang,
Backward stochastic partial differential equations with jumps and
application to optimal control of random jump fields, Stochastics,
77, No. 5 (2005), 381--399.
\bibitem{[Pe-93]} S. G. Peng, Backward stochastic differential
equations and applications to optimal control, Appl. Math. Optim.
27, No. 2 (1993), 125--144.
\bibitem{[Pe-Za07]} S. Peszat and J. Zabczyk, Stochastic
partial differential equations with Lévy noise. An evolution
equation approach, Encyclopedia of Mathematics and its Applications
113, Cambridge University Press, Cambridge, 2007.
\bibitem{[Pont61]} L. S. Pontryagin, Optimal regulation processes,
Amer. Math. Soc. Transl., 18, 2 (1961), 321--339.
\bibitem{[Roz90]} B. L. Rozovski\u\i,  Stochastic evolution systems.
Linear theory and applications to nonlinear filtering, Translated
from the Russian by A. Yarkho. Mathematics and its Applications
(Soviet Series), 35, Kluwer Academic Publishers Group, Dordrecht,
1990.
\bibitem{[Tang-Li]} S. Tang and X. Li, Mximum principle for optimal control
of distributed parameter stochastic systems with random jumps,
Differential equations, dynamical systems, and control science,
867–890, Lecture Notes in Pure and Appl. Math., 152, Dekker, New York, 1994.
\bibitem{[Tud89]} C. Tudor, Optimal control for semilinear stochastic
evolution equations, Appl. Math. Optim., 20, No. 3 (1989), 319--331.
\bibitem{[Y-Z]} J. Yong and X. Y. Zhou, Stochastic controls.
Hamiltonian systems and HJB equations, Springer-Verlag, New-York,
1999.
\bibitem{[Zhou93]} X. Y. Zhou, On the necessary conditions of optimal controls
for stochastic partial differential equations,
SIAM J. Control Optim., 31 (1993), No. 6, 1462--1478.

\end{thebibliography}
\end{document}